# Approximation of smooth functions on compact two-point homogeneous spaces

Gavin Brown[a], Feng Dai[b,*,1]

[a]*University of Sydney, NSW 2006, Australia*
[b]*Department of Mathematical and Statistical Sciences, Central Academic Building 632, University of Alberta, Edmonton, Alberta Canada T6G 2G1*

**Abstract**

Estimates of Kolmogorov $n$-widths $d_n(B_p^r, L^q)$ and linear $n$-widths $\delta_n(B_p^r, L^q)$, $(1 \leqslant q \leqslant \infty)$ of Sobolev's classes $B_p^r$, $(r>0, 1 \leqslant p \leqslant \infty)$ on compact two-point homogeneous spaces (CTPHS) are established. For part of $(p,q) \in [1, \infty] \times [1, \infty]$, sharp orders of $d_n(B_p^r, L^q)$ or $\delta_n(B_p^r, L^q)$ were obtained by Bordin et al. (J. Funct. Anal. 202(2) (2003) 307). In this paper, we obtain the sharp orders of $d_n(B_p^r, L^q)$ and $\delta_n(B_p^r, L^q)$ for all the remaining $(p,q)$. Our proof is based on positive cubature formulas and Marcinkiewicz–Zygmund-type inequalities on CTPHS.
© 2004 Elsevier Inc. All rights reserved.

## 1. Introduction

Given an integer $d \geqslant 3$, we denote by $M^{d-1}$ a compact two-point homogeneous space (CTPHS) of dimension $d-1$. According to Wang Hsien-Chung [W], the spaces of this type are the unit spheres $\mathbb{S}^{d-1}$, $d = 2, 3, \ldots$, the real projective spaces $P^{d-1}(R)$, $d = 3, 4, \ldots$, the complex projective spaces $P^{d-1}(C)$, $d = 5, 7, \ldots$, the quaternion projective spaces $P^{d-1}(H)$, $d = 9, 13, 17, \ldots$ and the Cayley projective plane $P^{d-1}(Cay)$, $d = 17$. These spaces are the compact symmetric spaces of rank one.

---

*∗ Corresponding author.*
*E-mail addresses:* vice-chancellor@vcc.usyd.edu.au (G. Brown), dfeng@math.ualberta.ca (F. Dai).
[1] Conducted this work as a student at the University of Sydney with support from the Australian Research Council.

We denote by $d\sigma(x)$ the Riemannian measure on $M^{d-1}$ normalized by $\int_{M^{d-1}} d\sigma(x) = 1$. Also we denote by $\|f\|_p$ the quantity $\left(\int_{M^{d-1}} |f(x)|^p \, d\sigma(x)\right)^{\frac{1}{p}}$ for $0 < p < \infty$ and $\|f\|_\infty$ the essential supremum of $f$ over $M^{d-1}$. Thus, for $0 < p \leqslant \infty$, $L^p \equiv L^p(M^{d-1})$ is a linear space equipped with the quasi-norm $\|\cdot\|_p$. Given $1 \leqslant p \leqslant \infty$ and $r > 0$, we denote by $B_p^r \equiv B_p^r(M^{d-1})$ the Sobolev class of order $r$ on $M^{d-1}$. (The precise definition of $B_p^r$ will be given in Section 2.)

For a given subset $K$ of a normed linear space $(X, \|\cdot\|)$, the Kolmogorov $n$-width $d_n(K, X)$ is defined by

$$d_n(K, X) = \inf_{L_n} \sup_{x \in K} \inf_{y \in L_n} \|x - y\|,$$

with the left-most infimum being taken over all $n$-dimensional linear subspaces $L_n$ of $X$, while the linear $n$-width $\delta_n(K, X)$ is defined by

$$\delta_n(K, X) = \inf_{T_n} \sup_{x \in K} \|x - T_n(x)\|,$$

with the infimum being taken over all linear continuous operators $T_n$ from $X$ to $X$ with $\dim\left(T_n(X)\right) \leqslant n$.

Our main goal in this paper is to determine the sharp asymptotic orders of $d_n(B_p^r, L^q)$ and $\delta_n(B_p^r, L^q)$ as $n \to \infty$. Before stating our results, we would like to sketch briefly some related known results. First, in the case when $M^{d-1} = \mathbb{S}^{d-1}$ (the unit sphere of the $d$-dimensional Euclidean space $\mathbb{R}^d$), the investigation of this problem began with the work of Kamzolov [Ka1,Ka2], who proved the exact orders of $d_n(B_p^r(\mathbb{S}^{d-1}), L^q(\mathbb{S}^{d-1}))$ for $1 \leqslant p = q \leqslant \infty$ or $1 \leqslant p \leqslant q \leqslant 2$; incorrect proofs of the orders of $d_n(B_p^r(\mathbb{S}^{d-1}), L^q(\mathbb{S}^{d-1}))$ and $\delta_n(B_p^r(\mathbb{S}^{d-1}), L^q(\mathbb{S}^{d-1}))$ were presented by Kushpel [Ku], as was pointed out in [BDS]; assuming the stronger smoothness condition $r \geqslant 2(d-1)^2$, Sun Yongsheng and the current authors obtained the exact orders of $d_n(B_p^r(\mathbb{S}^{d-1}), L^q(\mathbb{S}^{d-1}))$ for $1 \leqslant p \leqslant 2 \leqslant q \leqslant \infty$ or $2 \leqslant p \leqslant q \leqslant \infty$ in [BDS]. Second, for a general CTPHS $M^{d-1}$, the exact orders of $d_n$ and $\delta_n$ for $1 \leqslant p = q \leqslant \infty$, or $2 \leqslant q \leqslant p < \infty$, or $1 \leqslant p \leqslant q \leqslant 2$, and those of $\delta_n$ for $1 < q \leqslant p \leqslant 2$ or $2 \leqslant p \leqslant q \leqslant \infty$ were obtained in a recent paper [BKLT]. In this paper, we will determine the orders of $d_n(B_p^r(M^{d-1}), L^q(M^{d-1}))$ and $\delta_n(B_p^r(M^{d-1}), L^q(M^{d-1}))$ for all the remaining $(p, q)$. Moreover, in our work we allow the smoothness condition on $r$ to be relaxed to the expected one. Our main results (Theorem 1.1 below) will be in full analogy with those of Kashin [Kas], Maiorov [Ma] and Höllig [Ho] for the case of unit circle.

**Theorem 1.1.** *For* $1 \leqslant p, q \leqslant \infty$, $p' = \frac{p}{p-1}$, *we have*

$$d_n(B_p^r, L^q) \asymp \begin{cases} n^{-\frac{r}{d-1} + \frac{1}{p} - \frac{1}{2}} & \text{if } 1 \leqslant p \leqslant 2 \leqslant q \leqslant \infty \text{ and } r > \frac{d-1}{p}, \\ n^{-\frac{r}{d-1}} & \text{if } 2 \leqslant p \leqslant q \leqslant \infty \text{ and } r > \frac{d-1}{p}, \\ n^{-\frac{r}{d-1}} & \text{if } 1 \leqslant q \leqslant p \leqslant \infty \text{ and } r > 0 \end{cases}$$

*and*

$$\delta_n(B_p^r, L^q) \asymp \begin{cases} n^{-\frac{r}{d-1} + \frac{1}{p} - \frac{1}{2}} & \text{if } 1 \leqslant p \leqslant 2 \leqslant q \leqslant p' \text{ and } r > \frac{d-1}{p}, \\ n^{-\frac{r}{d-1} + \frac{1}{2} - \frac{1}{q}} & \text{if } 1 \leqslant p \leqslant 2 \leqslant p' \leqslant q \leqslant \infty \text{ and } r > (d-1)(1 - \frac{1}{q}), \\ n^{-\frac{r}{d-1}} & \text{if } 1 \leqslant q \leqslant p \leqslant \infty \text{ and } r > 0, \end{cases}$$

*where* '$\asymp$' *means that the ratio of both sides lies between two positive constants which are independent of* $n$.

We organize this paper as follows. Section 2 contains some basic notations and facts concerning harmonic analysis on a CTPHS $M^{d-1}$. In Section 3, we establish positive cubature formulas and Marcinkiewicz–Zygmund (MZ) inequalities on $M^{d-1}$, which give an extension of the results obtained in [MNW1]. Finally, based on the results obtained in Section 3, we prove Theorem 1.1 in Section 4.

Throughout the paper, the letter $C$ denotes a general positive constant depending only on the parameters indicated as subscripts, and the notation $A \asymp B$ means that there are two inessential positive constants $C_1$, $C_2$ such that $C_1 A \leqslant B \leqslant C_2 A$.

## 2. Harmonic analysis on $M^{d-1}$

This section is devoted to a brief description of some basic facts and notations concerning harmonic analysis on a CTPHS $M^{d-1}$. Most of the material in this section can be found in [BC, Section 7] and [He2, Chapter I, Section 4].

We denote by $d(\cdot, \cdot)$ the Riemannian metric on $M^{d-1}$ normalized so that all geodesics on $M^{d-1}$ have the same length $2\pi$. We denote by $B(x, r)$ the ball centered at $x \in M^{d-1}$ and having radius $r > 0$, i.e., $B(x, r) = \{y \in M^{d-1} : d(x, y) \leqslant r\}$, and $|E|$ the measure $\sigma(E)$ of a measurable subset $E \subset M^{d-1}$.

The space $M^{d-1}$ can be identified as a homogeneous space $G/K$, where $G$ is the maximal connected group of isometries of $M^{d-1}$ and $K = \{T \in G : Te = e\}$, $e$ is a fixed element in $M^{d-1}$. It is known that both the measure $d\sigma(x)$ and the metric $d(\cdot, \cdot)$ on $M^{d-1}$ are $G$-invariant and for any fixed $o \in M^{d-1}$ and any $f \in L(M^{d-1})$,

$$\int_{M^{d-1}} f(x) \, d\sigma(x) = \int_G f(go) \, dg,$$

where $dg$ is the usual Haar measure on $G$ normalized by $\int_G dg = 1$.

Given $t \in (0, \pi)$, let $\alpha(t)$ denote the total volume of a sphere of radius $t$ in $M$. It is well known (see [He2, p. 168]) that

$$\alpha(t) = C \left(\sin \frac{t}{2}\right)^a (\sin t)^b,$$

where $C > 0$ is chosen so that $\int_0^\pi \alpha(t)\, dt = 1$, $a$ and $b$ are determined as follows:

$$M^{d-1} = \mathbb{S}^{d-1}: \quad a = 0, \ b = d-2; \quad M^{d-1} = P^{d-1}(R): \quad a = d-2, \ b = 0,$$
$$M^{d-1} = P^{d-1}(C): \quad a = d-3, b = 1; \quad M^{d-1} = P^{d-1}(H): \quad a = d-5, b = 3,$$
$$M^{d-1} = P^{16}(Cay): \quad a = 8, b = 7. \tag{2.1}$$

Moreover, if $f(t)$ is integrable on $\big([0, \pi], \alpha(t)\, dt\big)$, then for any fixed $o \in M^{d-1}$,

$$\int_{M^{d-1}} f(d(x, o))\, d\sigma(x) = \int_0^\pi f(t)\alpha(t)\, dt. \tag{2.2}$$

Since $a + b = d - 2$, it follows that for any $o \in M^{d-1}$ and any $r \in (0, \pi)$,

$$|B(o, r)| = \int_0^r \alpha(t)\, dt \asymp r^{d-1} \tag{2.3}$$

and that for any $o \in M^{d-1}$, $0 < \gamma < \gamma + r \leqslant \pi$,

$$\left|\{x \in M^{d-1}: \gamma \leqslant d(x, o) \leqslant \gamma + r\}\right| = \int_\gamma^{\gamma+r} \alpha(t)\, dt \leqslant C_d (\gamma + r)^{d-2} r. \tag{2.4}$$

For convenience, for the rest of this section, we set $\varepsilon = 1$ or $2$ according as $M^{d-1} \neq P^{d-1}(R)$ or $M^{d-1} = P^{d-1}(R)$, and set

$$\alpha = \frac{d-3}{2}, \quad \beta = \frac{(d-2)(\varepsilon - 1) + b - 1}{2} \tag{2.5}$$

with $b$ being determined by (4.1).

Let $\triangle$ denote the Laplace–Beltrami operator on $M^{d-1}$. The spectrum of $\triangle$ is discrete, real and non-positive. We arrange it in decreasing order

$$0 = \lambda_0 > \lambda_1 > \lambda_2 > \cdots$$

and denote by $\mathcal{H}_k$ the eigenspace of $\triangle$ corresponding to the eigenvalue $\lambda_k$. It is well known that $\lambda_k = -\varepsilon k(\varepsilon k + \alpha + \beta + 1)$, the spaces $\mathcal{H}_k$, $(k \geqslant 0)$ are mutually orthogonal with respect to the inner product $\langle f, g \rangle := \int_{M^{d-1}} f(y)\overline{g(y)} \, d\sigma(y)$, and moreover $L^2(M) = \bigoplus_{k=0}^{\infty} \mathcal{H}_k$. For each integer $k \geqslant 0$ and each orthonormal basis of $\mathcal{H}_k$:

$$S_{k,1}(x), \ldots, S_{k,m_k}(x), \quad m_k = \dim \mathcal{H}_k,$$

we have the following addition formula (see [Gi]):

$$\sum_{j=1}^{m_k} S_{k,j}(x)\overline{S_{k,j}(y)} = c_{\varepsilon k} P_{\varepsilon k}^{(\alpha,\beta)}\left(\cos(\varepsilon^{-1} d(x,y))\right), \quad x, y \in M^{d-1}, \tag{2.6}$$

where

$$c_k = \frac{\Gamma(\beta+1)(2k+\alpha+\beta+1)\Gamma(k+\alpha+\beta+1)}{\Gamma(\alpha+\beta+2)\Gamma(k+\beta+1)} \tag{2.7}$$

and throughout, $P_k^{(\alpha_1,\beta_1)}$, $(\alpha_1, \beta_1 > -1)$ denotes the Jacobi polynomial normalized by

$$P_k^{(\alpha_1,\beta_1)}(1) = \frac{\Gamma(k+\alpha_1+1)}{\Gamma(k+1)\Gamma(\alpha_1+1)}. \tag{2.8}$$

(For the precise definition of $P_k^{(\alpha_1,\beta_1)}$ we refer to [Sz, p. 58].) The addition formula allows us to express the orthogonal projection $Y_k$ of $L^2(M^{d-1})$ onto $\mathcal{H}_k$ as the following convolution:

$$Y_k(f)(x) = c_{\varepsilon k} \int_{M^{d-1}} P_{\varepsilon k}^{(\alpha,\beta)}(\cos(\varepsilon^{-1} d(x,y))) f(y) \, d\sigma(y), \quad x \in M^{d-1}. \tag{2.9}$$

Evidently, from this last equality, the definition of $Y_k(f)$ can be extended to include all distributions $f$ on $M^{d-1}$.

For an integer $N \geqslant 0$, we put $\Pi_N \equiv \Pi_N(M^{d-1}) = \bigoplus_{k=0}^{N} \mathcal{H}_k$. The functions in $\Pi_N$ are called spherical polynomials of degree at most $N$. In the case $M^{d-1} = \mathbb{S}^{d-1}$, these functions coincide with the ordinary spherical polynomials (i.e. polynomials in $d$-variables restricted to $\mathbb{S}^{d-1}$). It is known that for any $f, g \in \Pi_N(M^{d-1})$, $fg \in \Pi_{2N}(M^{d-1})$. (For $\mathbb{S}^{d-1}$, see [WL, Chapter I]; for $P^{d-1}(R)$, this can be deduced directly from the case of $\mathbb{S}^{d-1}$; for $P^{d-1}(C)$, $P^{d-1}(H)$ and $P^{16}(Cay)$, this is a consequence of [Gr]). The addition formula implies that $\dim \mathcal{H}_k = c_{\varepsilon k} P_{\varepsilon k}^{(\alpha,\beta)}(1) \asymp k^{d-2}$ and hence $\dim \Pi_N = \sum_{k=0}^{N-1} \dim \mathcal{H}_k \asymp N^{d-1}$.

Given $r > 0$, we define the $r$th order Laplace–Beltrami operator $(-\Delta)^r$ on $M^{d-1}$ in a distributional sense by

$$Y_k\Big((-\Delta)^r(f)\Big) := \Big(\varepsilon k(\varepsilon k + \alpha + \beta + 1)\Big)^r Y_k(f), \quad k = 0, 1, \ldots,$$

where $f$ is a distribution on $M^{d-1}$. For $r > 0$ and $1 \leqslant p \leqslant \infty$, the Sobolev space $W_p^r \equiv W_p^r(M^{d-1})$ is defined by

$$W_p^r := \Big\{ f \in L^p(M^{d-1}) : \quad (-\Delta)^{r/2}(f) \in L^p(M^{d-1}) \Big\},$$

while the Sobolev class $B_p^r \equiv B_p^r(M^{d-1})$ is defined to be the set of all functions $f \in W_p^r$ such that $\|(-\Delta)^{r/2}(f)\|_p \leqslant 1$.

For the reminder of this section, we assume $M^{d-1} \neq P^{d-1}(R)$.

Given $f \in L^1(M^{d-1})$, we define the Cesàro means $\sigma_k^\delta(f)$, ($k \in \mathbb{Z}_+$) of order $\delta > -1$ by

$$\sigma_k^\delta(f) := \sum_{j=0}^{k} \frac{A_{k-j}^\delta}{A_k^\delta} Y_j(f),$$

where $A_j^\delta = \frac{\Gamma(j+\delta+1)}{\Gamma(\delta+1)\Gamma(j+1)}$. From [BC, Section 7] we know that for $\delta > \frac{d-2}{2}$ and $1 \leqslant p \leqslant \infty$,

$$\sup_{k \in \mathbb{Z}_+} \|\sigma_k^\delta(f)\|_p \leqslant C_\delta \|f\|_p. \tag{2.10}$$

Let $\eta$ be a $C^\infty$-function on $[0, \infty)$ supported in $[0, 2]$ and being equal to 1 on $[0, 1]$. For each integer $N \geqslant 1$, we define

$$V_{N,\eta}(f) := \sum_{k=0}^{2N} \eta\left(\frac{k}{N}\right) Y_k(f) \tag{2.11}$$

and

$$K_{N,\eta}(u) := \sum_{k=0}^{2N} \eta\left(\frac{k}{N}\right) c_k P_k^{(\alpha,\beta)}(u), \quad u \in [-1, 1], \tag{2.12}$$

where $c_k$ is defined by (2.7), $\alpha, \beta$ are determined by (2.5). From (2.9) and (2.10), it is easily seen that for $P \in \Pi_N$,

$$V_{N,\eta}(P)(x) = \int_{M^{d-1}} P(y) K_{N,\eta}(\cos d(x,y)) \, d\sigma(y) = P(x) \tag{2.13}$$

and moreover, for all $1 \leq p \leq \infty$ and $f \in L^p(M^{d-1})$,

$$\sup_N \|V_{N,\eta}(f)\|_p \leq C_{d,\eta} \|f\|_p. \tag{2.14}$$

We will keep the notations $\alpha, \beta, \alpha(t), \eta, K_{N,\eta}$ and $V_{N,\eta}$ for the rest of the paper.

## 3. Positive cubature formulas and MZ inequalities on $M^{d-1}$

Given a finite set $\Lambda$, we denote by $\#\Lambda$ the cardinality of $\Lambda$. For $r > 0$ and $a \geq 1$, we say a finite subset $\Lambda \subset M^{d-1}$ is an $(r,a)$-covering of $M^{d-1}$ if

$$M^{d-1} \subset \bigcup_{\omega \in \Lambda} B(\omega, r) \quad \text{and} \quad \max_{\omega \in \Lambda} \#\left(\Lambda \bigcap B(\omega, r)\right) \leq a.$$

By (2.3) and the definition, it is easily seen that every $(r,a)$-covering $\Lambda$ must satisfy the following condition

$$1 \leq \sum_{\omega \in \Lambda} \chi_{B(\omega,r)}(x) = \#\left(B(x,r) \bigcap \Lambda\right) \leq C_1 a \quad \text{for any } x \in M^{d-1}, \tag{3.1}$$

where $\chi_E$ denotes the characteristic function of a measurable set $E$, and $C_1$ is a constant depending only on $d$.

Our main results in this section can be stated as follows.

**Theorem 3.1.** *There exists a constant $\gamma > 0$ depending only on $d$ such that for any positive integer $n$ and any $(\frac{\delta}{n}, a)$-covering of $M^{d-1}$ satisfying $0 < \delta < a^{-1}\gamma$, there exists a set of numbers $0 \leq \lambda_\omega \leq C_d n^{-(d-1)}$, $(\omega \in \Lambda)$ such that*

$$\int_{M^{d-1}} f(y) \, d\sigma(y) = \sum_{\omega \in \Lambda} \lambda_\omega f(\omega) \quad \text{for any } f \in \Pi_{4n} \tag{3.2}$$

*and moreover, for any $0 < p \leq \infty$, $f \in \Pi_n$ and $0 \leq t \leq \min\{p, 1\}$,*

$$\|f\|_p \asymp \begin{cases} \left(\frac{1}{n^{d-1}} \sum_{\omega \in \Lambda} \left(n^{d-1} \lambda_\omega\right)^t |f(\omega)|^p\right)^{\frac{1}{p}} & \text{if } 0 < p < \infty, \\ \max_{\omega \in \Lambda} \left(n^{d-1} \lambda_\omega\right)^t |f(\omega)| & \text{if } p = \infty, \end{cases} \tag{3.3}$$

*where the constants of equivalence depend only on d, a, $\delta$ and p when p is small, and we employ the slight abuse of notation that $0^0 = 1$.*

An equality like (3.2) with nonnegative weights is generally called a positive cubature formula of degree $4n$, while a relation like (3.3) is called an MZ-type inequality.

Of special interest is the case when $M^{d-1} = \mathbb{S}^{d-1}$. In this case, an important breakthrough has been achieved recently by Mhaskar et al. [MNW1,MNW2], who obtained MZ inequalities and positive cubature formulas based on function values at scattered sites on the unit sphere $\mathbb{S}^{d-1}$. Although our work in this section was greatly influenced by the remarkable paper [MNW1], we believe it is of independent interest because of the following three reasons. First, the MZ inequalities in our theorem apply to the full range of $p$, i.e., $0 < p \leq \infty$, while those in [MNW1] were established for the case $1 \leq p \leq \infty$. Second, we allow more choices of the weights in our MZ inequalities. Indeed, these weights can be taken to be equal or to be the same as those in the cubature formulas when $p \geq 1$. This will play a very important role in our proof of Theorem 1.1 in Section 4. (A referee kindly informed us that the techniques in [MNW1] can also apply to obtain 'one-weight' MZ inequalities on $\mathbb{S}^{d-1}$.) Third, the proof in [MNW1] is based on an ingenious decomposition of $\mathbb{S}^{d-1}$ into a finite collection of comparable spherical simplices and involves some doubling-weight tricks, while our proof here is different and less complicated. More importantly, it works for all CTPHS.

For the proof of Theorem 3.1, without loss of generality, we may assume $M^{d-1} \neq P^{d-1}(R)$. Indeed, since the real projective space $P^{d-1}(R)$ is the quotient of $\mathbb{S}^{d-1}$ by the equivalence relation which identifies $x$ with $-x$, functions on $P^{d-1}(R)$ can be thought of as even functions on $\mathbb{S}^{d-1}$, in which case $\mathcal{H}_k(P^{d-1}(R)) = \mathcal{H}_{2k}(\mathbb{S}^{d-1})$ for all $k \geq 0$. Thus, Theorem 3.1 for the case $M^{d-1} = P^{d-1}(R)$ follows directly from the case when $M^{d-1} = \mathbb{S}^{d-1}$.

We remind the reader that throughout the rest of this paper the notations $\alpha$, $\beta$, $\alpha(t)$, $\eta$, $V_{N,\eta}$ and $K_{N,\eta}$ will be the same as those in Section 2.

The proof of Theorem 3.1 depends on the following two lemmas.

**Lemma 3.2.** *Suppose that $\delta \in (0, \pi)$, $a \geq 1$, n is a positive integer and $\Lambda$ is a $(\frac{\delta}{n}, a)$-covering of $M^{d-1}$. Then for all $1 \leq p < \infty$ and all $f \in \Pi_{4n}$,*

$$\left( \sum_{\omega \in \Lambda} \left| B\left(\omega, \frac{\delta}{n}\right) \right| \max_{x \in B(\omega, \frac{\delta}{n})} |f(x) - f(\omega)|^p \right)^{\frac{1}{p}} \leq C_2 a^{\frac{1}{p}} \delta \|f\|_p, \tag{3.4}$$

*where the constant $C_2$ depends only on d.*

We point out that for the unit sphere $\mathbb{S}^{d-1}$ an $L^1$-inequality similar to (3.4) was previously obtained in [MNW1].

**Lemma 3.3.** *For any $\ell \geqslant 1$ and $\theta \in (0, \pi)$,*

$$|K_{N,\eta}^{(i)}(\cos\theta)| \leqslant C_{\ell,\eta,i} N^{d-1+2i} \min\{1, (N\theta)^{-\ell}\}, \quad i = 0, 1, 2, \ldots, \qquad (3.5)$$

where $K_{N,\eta}^{(0)}(u) = K_{N,\eta}(u)$ and $K_{N,\eta}^{(i)}(u) = \left(\frac{d}{du}\right)^i \{K_{N,\eta}(u)\}$ for $i \geqslant 1$.

We note that for the Cesàro kernel $\sigma_N^\delta(\cos\theta) = \sum_{k=0}^N \frac{A_{N-k}^\delta}{A_N^\delta} c_k P_k^{(\alpha,\beta)}(\cos\theta)$ with $A_k^\delta = \frac{\Gamma(k+\delta+1)}{\Gamma(k+1)\Gamma(\delta+1)}$ and $c_k$ defined by (2.7), it is well known that (see [BC, Theorem 2.1]) for $\delta > \alpha + \beta + 2$, $|\sigma_N^\delta(\cos\theta)| \leqslant C_\delta N^{d-1} \min\{1, (N\theta)^{-d}\}$ and moreover, (see [BC, p. 234, Remark 1]) the order $(N\theta)^{-d}$ on the right-hand side of this last estimate is best possible. The significant point of our Lemma 3.3 is that the number $\ell$ can be chosen as large as we like, which will play a very important role when we prove Theorem 3.1 for $0 < p < 1$.

For the moment, we take Lemmas 3.2 and 3.3 for granted and proceed with the proof of Theorem 3.1. The proofs of these two lemmas will be given in Sections 3.3 and 3.2, respectively.

### 3.1. Proof of Theorem 3.1

Throughout the proof in this subsection, $C_1$ and $C_2$ will denote the constants appearing in (3.1) and (3.4), respectively. Let $\gamma = (2C_2)^{-1}$ and suppose $\Lambda$ is a $(\frac{\delta}{n}, a)$-covering of $M^{d-1}$ with $a \geqslant 1$ and $0 < \delta < \gamma a^{-1}$. We will prove that $\Lambda$ has the desired properties of Theorem 3.1.

First, we show that there exists a set of numbers $0 \leqslant \lambda_\omega \leqslant C_d n^{-(d-1)}$, $(\omega \in \Lambda)$ such that the cubature formula (3.2) holds for all $f \in \Pi_{4n}$. The proof is based on Lemma 3.2 while the idea is from Mhaskar et al. [MNW1]. In fact, using Lemma 3.2 and (3.1) with $r = \frac{\delta}{n}$, it is easy to verify that for $1 \leqslant p < \infty$ and $f \in \Pi_{4n}$,

$$\frac{1}{2}\left(1 - (2C_2\delta)^p a\right)^{\frac{1}{p}} \|f\|_p \leqslant \left(\sum_{\omega \in \Lambda} \left|B\left(\omega, \frac{\delta}{n}\right)\right| |f(\omega)|^p\right)^{\frac{1}{p}}$$

$$\leqslant 2\left((C_2\delta)^p a + C_1\right)^{\frac{1}{p}} \|f\|_p \qquad (3.6)$$

and moreover, for $f \in \Pi_{4n}$ with $\min_{\omega \in \Lambda} f(\omega) \geqslant 0$,

$$\int_{M^{d-1}} f(x)\, d\sigma(x) \geqslant \|f\|_1 - 2 \sum_{\omega \in \Lambda} \left|B\left(\omega, \frac{\delta}{n}\right)\right| \max_{x \in B(\omega, \frac{\delta}{n})} |f(x) - f(\omega)|$$

$$\geqslant (1 - 2C_2 a\delta)\|f\|_1. \qquad (3.7)$$

Now following the proof of Theorem 4.1 of [MNW1], using estimates (3.6) and (3.7), and by the "Krein–Rutman theorem" (see [Hol, p. 20]), we conclude that there must exist a set of nonnegative numbers $\lambda_\omega$, ($\omega \in \Lambda$) for which the cubature formula (3.2) holds for all $f \in \Pi_{4n}$. (We refer to [MNW1, Section 4] for the detailed proof.) For the proof of the inequality $\lambda_\omega \leqslant C_d n^{-(d-1)}$, we let $P_n(u) = \frac{P_n^{(\alpha,\beta)}(u)}{P_n^{(\alpha,\beta)}(1)}$ with $\alpha, \beta$ determined by (2.5), and apply the cubature formula (3.2) with $f(\cdot) = \left(P_n(\cos d(\cdot, \omega))\right)^2$ to obtain

$$\lambda_\omega \leqslant \sum_{\omega' \in \Lambda} \lambda_{\omega'} \left(P_n(\cos d(\omega', \omega))\right)^2 = \int_{M^{d-1}} \left(P_n(\cos d(x, \omega))\right)^2 d\sigma(x)$$
$$= \int_0^\pi \left(P_n(\cos \theta)\right)^2 \alpha(\theta) \, d\theta,$$

which is dominated by the desired bound $C_d n^{-(d-1)}$, on account of (2.8) and the following well-known estimates on Jacobi polynomials [Sz, (7.32.6)], p. 167]: for $k \geqslant 1$ and $\alpha_1, \beta_1 > -\frac{1}{2}$,

$$|P_k^{(\alpha_1, \beta_1)}(\cos \theta)| \leqslant C_{\alpha_1, \beta_1} \begin{cases} \min\{k^{\alpha_1}, k^{-\frac{1}{2}} \theta^{-\alpha_1 - \frac{1}{2}}\} & \text{if } 0 \leqslant \theta \leqslant \frac{\pi}{2}, \\ \min\{k^{\beta_1}, k^{-\frac{1}{2}}(\pi - \theta)^{-\beta_1 - \frac{1}{2}}\} & \text{if } \frac{\pi}{2} \leqslant \theta \leqslant \pi. \end{cases} \quad (3.8)$$

Second, we prove (3.3) for $1 \leqslant p \leqslant \infty$, $t \in [0, 1]$ and $f \in \Pi_n$. By (3.6) and the inequality

$$0 \leqslant \lambda_\omega \leqslant C_d n^{-(d-1)}, \quad (3.9)$$

it will suffice to prove that for $f \in \Pi_n$

$$\|f\|_p \leqslant C_{d,a,\delta} \begin{cases} \left(\sum_{\omega \in \Lambda} \lambda_\omega |f(\omega)|^p\right)^{\frac{1}{p}} & \text{if } 1 \leqslant p < \infty, \\ \max_{\omega \in \Lambda}\left((\lambda_\omega n^{d-1})|f(\omega)|\right) & \text{if } p = \infty. \end{cases} \quad (3.10)$$

This can be obtained by duality. In fact, by the cubature formula (3.2) and the mutual orthogonity of the spaces $\mathcal{H}_k$, $k \in \mathbb{Z}_+$, it follows that for $f \in \Pi_n$ and $g \in L^{p'}$ with $\frac{1}{p} + \frac{1}{p'} = 1$,

$$\int_{M^{d-1}} f(x) g(x) \, d\sigma(x) = \int_{M^{d-1}} f(x) V_{n,\eta}(g)(x) \, d\sigma(x) = \sum_{\omega \in \Lambda} \lambda_\omega f(\omega) V_{n,\eta}(g)(\omega).$$

Hence, by Hölder's inequality, (3.9), (2.14), (2.3) with $r = \frac{\delta}{n}$, and (3.6) applied to $p'$ and $V_{n,\eta}(g)$, we obtain

$$\left| \int_{M^{d-1}} f(x) g(x) \, d\sigma(x) \right| \leqslant C_d \left( \sum_{\omega \in \Lambda} \lambda_\omega |f(\omega)|^p \right)^{\frac{1}{p}} \left( \frac{1}{n^{d-1}} \sum_{\omega \in \Lambda} |V_{n,\eta}(g)(\omega)|^{p'} \right)^{\frac{1}{p'}}$$

$$\leqslant C_{d,a,\delta} \|V_{n,\eta}(g)\|_{p'} \left( \sum_{\omega \in \Lambda} \lambda_\omega |f(\omega)|^p \right)^{\frac{1}{p}}$$

$$\leqslant C_{d,a,\delta} \|g\|_{p'} \left( \sum_{\omega \in \Lambda} \lambda_\omega |f(\omega)|^p \right)^{\frac{1}{p}},$$

with slight changes when $p = 1$ or $\infty$. Taking the supremum over all $g \in L^{p'}$ with $\|g\|_{p'} = 1$, we obtain (3.10) and hence complete the proof of (3.3) for $1 \leqslant p \leqslant \infty$.

Finally, we show (3.3) for $0 < p < 1$, $t \in [0, p]$ and $f \in \Pi_n$. By (2.13) and the cubature formula (3.2), we have, for $f \in \Pi_n$,

$$f(x) = \sum_{\omega \in \Lambda} \lambda_\omega K_{n,\eta}(\cos d(x, \omega)) f(\omega)$$

and hence for $0 < p < 1$,

$$\int_{M^{d-1}} |f(x)|^p \, d\sigma(x) \leqslant \sum_{\omega \in \Lambda} (\lambda_\omega)^p |f(\omega)|^p \int_{M^{d-1}} |K_{n,\eta}(\cos d(x, \omega))|^p \, d\sigma(x), \quad (3.11)$$

which, by (2.2) and Lemma 3.3 applied to $\ell > \frac{d-1}{p}$ and $i = 0$, is dominated by

$$C_{d,p} \sum_{\omega \in \Lambda} (\lambda_\omega n^{d-1})^p n^{-(d-1)} |f(\omega)|^p.$$

The desired upper estimate

$$\|f\|_p \leqslant C_{d,p} \left( \sum_{\omega \in \Lambda} (\lambda_\omega n^{d-1})^t n^{-(d-1)} |f(\omega)|^p \right)^{\frac{1}{p}}$$

for $t \in [0, p]$ and $f \in \Pi_n$ then follows since $0 \leqslant \lambda_\omega n^{d-1} \leqslant C_d$.

For the proof of the lower estimate

$$\left(\sum_{\omega\in\Lambda}(\lambda_\omega n^{d-1})^t n^{-(d-1)}|f(\omega)|^p\right)^{\frac{1}{p}}\leqslant C_{p,d,a,\delta}\|f\|_p \quad \text{for } t\in[0,p] \text{ and } f\in\Pi_n,$$

by (3.9), it will suffice to prove the case when $t = 0$. To this end, let $G$ denote the maximal connected group of isometries of $M^{d-1}$. Since the Riemannian measure $d\sigma(x)$ is $G$-invariant, it follows that for any $f \in \Pi_n$ and $g \in G$,

$$f(\omega)=\int_{M^{d-1}}f(gy)K_{n,\eta}(\cos d(\omega,gy))\,d\sigma(y)=\sum_{\omega'\in\Lambda}\lambda_{\omega'}f(g\omega')K_{n,\eta}(\cos d(\omega,g\omega')),$$

which together with the inequality $0\leqslant \lambda_{\omega'}\leqslant C_d n^{-(d-1)}$, $\omega' \in \Lambda$ implies

$$\frac{1}{n^{d-1}}\sum_{\omega\in\Lambda}|f(\omega)|^p$$

$$\leqslant C_{p,d}\inf_{g\in G}\left[\frac{1}{n^{(d-1)p}}\sum_{\omega'\in\Lambda}|f(g\omega')|^p\left(\frac{1}{n^{d-1}}\sum_{\omega\in\Lambda}|K_{n,\eta}(\cos d(\omega,g\omega'))|^p\right)\right]. \tag{3.12}$$

It will be shown that

$$\max_{y\in M^{d-1}}\left[\frac{1}{n^{d-1}}\sum_{\omega\in\Lambda}|K_{n,\eta}(\cos d(\omega,y))|^p\right]\leqslant C_{p,d}a\delta^{-(d-1)}n^{(d-1)(p-1)}, \tag{3.13}$$

which combined with (3.12) will imply the desired lower estimate:

$$\frac{1}{n^{d-1}}\sum_{\omega\in\Lambda}|f(\omega)|^p\leqslant C_{p,d,a,\delta}\frac{1}{n^{d-1}}\sum_{\omega'\in\Lambda}\int_G|f(g\omega')|^p\,dg\leqslant C_{p,d,a,\delta}\|f\|_p^p.$$

For the proof of (3.13), we set $\Lambda_0(y) = B(y,\frac{\pi}{n})\bigcap\Lambda$,

$$\Lambda_k(y)=\left\{\omega\in\Lambda:\ \frac{k\pi}{n}\leqslant d(\omega,y)\leqslant\frac{(k+1)\pi}{n}\right\},\quad 1\leqslant k\leqslant n-1$$

and then invoke Lemma 3.3 with $i = 0$ and $\ell = [\frac{d-1}{p}]+1$ to obtain

$$\frac{1}{n^{d-1}}\sum_{\omega\in\Lambda}|K_{n,\eta}(\cos d(\omega,y))|^p$$

$$\leqslant C_{p,d} \sum_{k=0}^{n-1} n^{(d-1)(p-1)} \sum_{\omega \in \Lambda_k(y)} (k+1)^{-p\ell}$$

$$= C_{p,d} \sum_{k=0}^{n-1} n^{(d-1)(p-1)} (k+1)^{-p\ell} \#\big(\Lambda_k(y)\big). \quad (3.14)$$

By the standard volumetric method, it is easily seen from (2.3) and (2.4) that $B(y, \frac{\pi}{n})$ can be covered by at most $C_d \delta^{-(d-1)}$ balls of radius $\frac{\delta}{2n}$ and that each set $\{z \in M^{d-1} : \frac{k\pi}{n} \leqslant d(z, y) \leqslant \frac{(k+1)\pi}{n}\}$ can be covered by at most $C_d k^{d-2} \delta^{-(d-1)}$ balls of radius $\frac{\delta}{2n}$. Since the definition of $(\frac{\delta}{n}, a)$-covering implies $\max_{x \in M^{d-1}} \big[\#\big(\Lambda \cap B(x, \frac{\delta}{2n})\big)\big] \leqslant a$, it follows that

$$\#\big(\Gamma_k(y)\big) \leqslant C_d \delta^{-(d-1)} a(k+1)^{d-2}, \quad 0 \leqslant k \leqslant n-1,$$

which combined with (3.14) implies (3.13) and therefore completes the proof of Theorem 3.1.  □

### 3.2. Proof of Lemma 3.3

Since (see [Sz, p. 63, (4.21.7)])

$$\frac{d}{dt}\{P_k^{(\alpha_1, \beta_1)}(t)\} = \frac{k + \alpha_1 + \beta_1 + 1}{2} P_{k-1}^{(\alpha_1+1, \beta_1+1)}(t) \quad \text{for } \alpha_1, \beta_1 > -1,$$

by the definition (2.12) it follows that for $i \in \mathbb{Z}_+$,

$$K_{N,\eta}^{(i)}(t) = C(\alpha, \beta, i) \sum_{k=0}^{2N-i} \eta\left(\frac{k+i}{N}\right)$$

$$\times \frac{(2k + \alpha + \beta + 2i + 1)\Gamma(k + \alpha + \beta + 2i + 1)}{\Gamma(k + i + \beta + 1)} P_k^{(\alpha+i, \beta+i)}(t).$$

(Without loss of generality, we may assume $N > i + \ell + 1$.) This together with the estimate (3.8) implies that for $\theta \in [0, N^{-1}]$

$$|K_{N,\eta}^{(i)}(\cos\theta)| \leqslant C_{d,\eta,i} \sum_{k=0}^{2N-i} (k+1)^{\alpha+i+1}(k+1)^{\alpha+i} \leqslant C_{d,\eta,i} N^{2\alpha+2i+2} = C_{d,\eta,i} N^{d-1+2i}.$$

Thus, it remains to prove that for any $\ell \geqslant 1$ and $\theta \in [N^{-1}, \pi]$,

$$|K_{N,\eta}^{(i)}(\cos\theta)| \leqslant C_{d,\eta,i,\ell} N^{d-1+2i}(N\theta)^{-\ell}. \quad (3.15)$$

To this end, we define a sequence $\{a_{N,j}(\cdot)\}_{j=0}^{\infty}$ of $C^{\infty}$-functions on $[0, \infty)$ by

$$a_{N,0}(u) = (2u + \alpha + \beta + 2i + 1)\eta\left(\frac{u+i}{N}\right),$$

$$a_{N,j+1}(u) = \frac{a_{N,j}(u)}{2u + \alpha + \beta + 2i + j + 1} - \frac{a_{N,j}(u+1)}{2u + \alpha + \beta + 2i + j + 3}, \quad j \geq 0.$$

Then, since for any integer $j \geq 0$, (see [Sz, (4.5.3), p. 71]),

$$\sum_{n=0}^{k} \frac{(2n + \alpha + \beta + 2i + j + 1)\Gamma(n + \alpha + \beta + 2i + j + 1)}{\Gamma(n + i + \beta + 1)} P_n^{(\alpha+i+j,\beta+i)}(t)$$

$$= \frac{\Gamma(k + \alpha + \beta + 2i + j + 2)}{\Gamma(k + \beta + i + 1)} P_k^{(\alpha+i+j+1,\beta+i)}(t),$$

it follows by summation by parts finite times that for any integer $j \geq 0$,

$$K_{N,\eta}^{(i)}(t) = C(\alpha, \beta, i) \sum_{k=0}^{\infty} a_{N,j}(k) \frac{\Gamma(k + \alpha + \beta + 2i + j + 1)}{\Gamma(k + \beta + i + 1)} P_k^{(\alpha+i+j,\beta+i)}(t). \quad (3.16)$$

Also, since $\operatorname{supp} \eta' \subset [1, 2]$, we can prove by induction that for $j \geq 1$

$$\operatorname{supp} a_{N,j}(\cdot) \subset [N - i - j, 2N - i], \quad (3.17)$$

$$\left|\left(\frac{d}{du}\right)^m a_{N,j}(u)\right| \leq C_{m,j,i,\varphi} N^{-(m+2j-1)}, \quad m = 0, 1, \ldots. \quad (3.18)$$

Therefore, using (3.16)–(3.18) with $m = 0$, taking into account (3.8), we obtain that for $\theta \in [N^{-1}, \pi/2]$ and $j \geq 1$,

$$|K_{N,\eta}^{(i)}(\cos\theta)| \leq C_{d,\eta,i,j} \sum_{k=N-i-j}^{2N-i} N^{-2j+1} k^{\alpha+i+j-\frac{1}{2}} \theta^{-(\alpha+i+j+\frac{1}{2})}$$

$$\leq C_{d,\eta,i,j} N^{d-1+2i} (N\theta)^{-(\frac{d-2}{2}+i+j)}$$

and for $\theta \in [\frac{\pi}{2}, \pi]$ and $j \geq 1$

$$|K_{N,\eta}^{(i)}(\cos\theta)| \leq C_{d,\eta,i,j} \sum_{k=N-i-j}^{2N-i} N^{-2j+1} k^{\alpha+i+j} k^{\beta+i} \leq C_{d,\eta,i,j} N^{\frac{d+1}{2}+\beta+2i-j}.$$

The desired estimate (3.15) then follows by taking $j = [\ell + \beta - \frac{d-3}{2}] + 1$. This completes the proof. $\square$

### 3.3. Proof of Lemma 3.2

For simplicity, we set, for $x, z \in M^{d-1}$ and $\delta > 0$,

$$G_{n,\delta}(x, z) := \max_{y \in B(x, \frac{\delta}{n})} |K_{4n,\eta}(\cos d(x, z)) - K_{4n,\eta}(\cos d(y, z))|. \tag{3.19}$$

Then using the triangle inequalities for the metric $d(\cdot, \cdot)$, and invoking Lemma 3.3 with $i = 0, 1$ and $\ell = d + 3$, we obtain that for $\delta \in (0, \pi)$ and $x, z \in M^{d-1}$,

$$G_{n,\delta}(x, z) \leqslant C_d \begin{cases} n^{d-1} & \text{if } 0 \leqslant d(x, z) \leqslant \frac{4\delta}{n}, \\ \delta n^{d-1} \min\{1, (nd(x, z))^{-(d+1)}\} & \text{if } \frac{4\delta}{n} \leqslant d(x, z) \leqslant \pi, \end{cases} \tag{3.20}$$

where $C_d$ is independent of $\delta$, $x$, $z$ and $n$.

We note that for $f \in \Pi_{4n}$ and $x, y \in M^{d-1}$,

$$f(x) - f(y) = \int_{M^{d-1}} f(z) \Big( K_{4n,\eta}(\cos d(x, z)) - K_{4n,\eta}(\cos d(y, z)) \Big) d\sigma(z).$$

Hence, using Hölder's inequality, (3.19) and (3.20), we obtain that for $1 \leqslant p < \infty$, $f \in \Pi_{4n}$, $\delta \in (0, \pi)$ and $x \in M^{d-1}$,

$$\max_{y \in B(x, \frac{\delta}{n})} |f(x) - f(y)|^p$$

$$\leqslant \left( \int_{M^{d-1}} |f(z)|^p G_{n,\delta}(x, z) d\sigma(z) \right) \left( \int_{M^{d-1}} G_{n,\delta}(x, z) d\sigma(z) \right)^{p-1}$$

$$\leqslant (C_d)^{p-1} \delta^{p-1} \int_{M^{d-1}} |f(z)|^p G_{n,\delta}(x, z) d\sigma(z).$$

It then follows that for $f \in \Pi_{4n}$ and $1 \leqslant p < \infty$,

$$\sum_{\omega \in \Lambda} \left| B\left(\omega, \frac{\delta}{n}\right) \right| \max_{y \in B(\omega, \frac{\delta}{n})} |f(\omega) - f(y)|^p$$

$$\leqslant (C_d)^{p-1} \delta^{p-1} \int_{M^{d-1}} |f(z)|^p \left( \sum_{\omega \in \Lambda} \left| B\left(\omega, \frac{\delta}{n}\right) \right| G_{n,\delta}(\omega, z) \right) d\sigma(z).$$

Therefore, it will suffice to prove

$$\max_{z \in M^{d-1}} \left( \sum_{\omega \in \Lambda} \left| B\left(\omega, \frac{\delta}{n}\right) \right| G_{n,\delta}(\omega, z) \right) \leqslant C_d \delta a. \tag{3.21}$$

To this end, we set, for $z \in M^{d-1}$,

$$\Gamma_{-1}(z) = \left\{\omega \in \Lambda : \ d(\omega, z) \leqslant \frac{4\delta}{n}\right\}, \quad \Gamma_0(z) := \left\{\omega \in \Lambda : \ \frac{4\delta}{n} \leqslant d(\omega, z) \leqslant \frac{\pi}{n}\right\}$$

and

$$\Gamma_k(z) := \left\{\omega \in \Lambda : \ \frac{k\pi}{n} \leqslant d(\omega, z) \leqslant \frac{(k+1)\pi}{n}\right\}, \quad 1 \leqslant k \leqslant n-1.$$

Then following the proof of (3.13), and using estimate (3.20), we have

$$\sum_{\omega \in \Lambda} \left|B\left(\omega, \frac{\delta}{n}\right)\right| G_{n,\delta}(\omega, z) \leqslant C_d \sum_{\omega \in \Gamma_{-1}(z)} n^{d-1} \left(\frac{\delta}{n}\right)^{d-1}$$

$$+ C_d \delta n^{d-1} \sum_{k=0}^{n-1} \sum_{\omega \in \Gamma_k(z)} \left(\frac{\delta}{n}\right)^{d-1} (k+1)^{-(d+1)}$$

$$= C_d \delta^{d-1} \#\big(\Gamma_{-1}(z)\big) + C_d \delta^d \sum_{k=0}^{n-1} (k+1)^{-(d+1)} \#\big(\Gamma_k(z)\big)$$

$$\leqslant C_d \delta^{d-1} a + C_d \delta^d \delta^{-(d-1)} a \sum_{k=1}^{n} k^{d-2} k^{-(d+1)} \leqslant C_d \delta a,$$

which gives (3.21) and thus completes the proof. □

## 4. Proof of Theorem 1.1

For $x = (x_1, \ldots, x_m) \in \mathbb{R}^m$ we define as usual $\|x\|_{\ell_p^m} = \left(\sum_{i=1}^m |x_i|^p\right)^{\frac{1}{p}}$ for $1 \leqslant p < \infty$ and $\|x\|_{\ell_\infty^m} = \max_{1 \leqslant i \leqslant m} |x_i|$. We denote by $\ell_p^m$ the set of vectors $x \in \mathbb{R}^m$ endowed with the norm $\|\cdot\|_{\ell_p^m}$ and $b_p^m$ the unit ball of $\ell_p^m$. Given $1 \leqslant p \leqslant \infty$ and an integer $N \geqslant 0$, we denote by $\mathcal{B}_N^p \equiv \mathcal{B}_N^p(M^{d-1})$ the class of all functions $f \in \Pi_N(M^{d-1})$ such that $\|f\|_p \leqslant 1$.

As Theorem 1.1 for the case $M^{d-1} = P^{d-1}(R)$ follows directly from the case when $M^{d-1} = \mathbb{S}^{d-1}$, we may assume $M^{d-1} \neq P^{d-1}(R)$ throughout this section.

The proof of Theorem 1.1 will rely on the following

**Lemma 4.1.** *Let $S_n$ denote either of the symbols $d_n$ or $\delta_n$. Then for $1 \leqslant p, q \leqslant \infty$ and $1 \leqslant n \leqslant \dim \Pi_N$, we have*

$$S_n(\mathcal{B}_N^p, L^q) \leqslant C_d N^{(d-1)(\frac{1}{p} - \frac{1}{q})} S_n(b_p^{m_N}, \ell_q^{m_N}),$$

*where $m_N \asymp \dim \Pi_N \asymp N^{d-1}$.*

**Proof.** Let $\gamma > 0$ be the same as in Theorem 3.1 and let $\Lambda = \{t_1, \ldots, t_m\}$ be a subset of $M^{d-1}$ such that $\min_{i \neq j} d(t_i, t_j) > \frac{\gamma}{10N}$ and $\max_{x \in M^{d-1}} \min_{1 \leq i \leq m} d(x, t_i) \leq \frac{\gamma}{10N}$. It is easily seen that $m = \#\Lambda \asymp N^{d-1}$ and $\Lambda$ is a $(\frac{\gamma}{10N}, 1)$-covering of $M^{d-1}$. Hence, by Theorem 3.1, there exists a set of numbers $0 \leq w_j \leq C_d N^{-(d-1)}$ such that for all $f \in \Pi_{8N}(M^{d-1})$, we have

$$\int_{M^{d-1}} f(y) \, d\sigma(y) = \sum_{j=1}^{m} w_j f(t_j) \tag{4.1}$$

and moreover, for all $1 \leq p \leq \infty$ and $f \in \Pi_{8N}(M^{d-1})$,

$$\|f\|_p \asymp N^{-\frac{d-1}{p}} \|U_N(f)\|_{\ell_p^m} \asymp \|U_N(f)\|_{\ell_{p,w}^m}, \tag{4.2}$$

where $U_N : \Pi_{8N} \to \mathbb{R}^m$ is defined by

$$U_N(f) = (f(t_1), \ldots, f(t_m)) \tag{4.3}$$

and for $u = (u_1, \ldots, u_m) \in \mathbb{R}^m$,

$$\|u\|_{\ell_{p,w}^m} = \begin{cases} \left(\sum_{j=1}^{m} |u_j|^p w_j\right)^{\frac{1}{p}} & \text{if } p < \infty, \\ \max_{1 \leq j \leq m} |u_j| & \text{if } p = \infty. \end{cases}$$

Next, we define

$$T(u)(\cdot) = \sum_{j=1}^{m} w_j u_j K_{N,\eta}\left(\cos(d(\cdot, t_j))\right), \quad u = (u_1, \ldots, u_m). \tag{4.4}$$

Since for a fixed $o \in M^{d-1}$ each $P_k^{(\alpha,\beta)}(\cos d(o, \cdot))$ is in $\mathcal{H}_k$, it follows that for each $o \in M^{d-1}$, $K_{N,\eta}\left(\cos d(o, \cdot)\right) \in \Pi_{2N}$, and hence $T$ is an operator from $\mathbb{R}^m$ to $\Pi_{2N}$. We will show that for $1 \leq q \leq \infty$,

$$\|T(u)(\cdot)\|_q \leq C_d \|u\|_{\ell_{q,w}^m}. \tag{4.5}$$

In fact, for $q = 1$, (4.5) follows directly from the inequality

$$\int_{M^{d-1}} |K_{N,\eta}(\cos d(x, o))| \, d\sigma(x) = \int_0^{\pi} |K_{N,\eta}(\cos \theta)| \alpha(\theta) \, d\theta \leq C_{d,\eta},$$

which can be easily deduced from Lemma 3.3; for $q = \infty$, by the second equivalence in (4.2) applied to $f(\cdot) = K_{N,\eta}(\cos d(x, \cdot))$,

$$\|T(u)\|_\infty \leqslant \|u\|_\infty \max_{x \in M^{d-1}} \sum_{j=1}^m w_j \left| K_{N,\eta}(\cos d(x, t_j)) \right|$$

$$\leqslant C_d \|u\|_\infty \max_x \int_{M^{d-1}} \left| K_{N,\eta}(\cos d(x, y)) \right| d\sigma(y)$$

$$= C_d \|u\|_\infty \int_0^\pi |K_{N,\eta}(\cos\theta)| \alpha(\theta) \, d\theta \leqslant C_d \|u\|_\infty$$

and for $1 < q < \infty$, (4.5) follows by the Riesz–Thorin theorem.

Finally, we use the cubature formula (4.1) to obtain that for $f \in \Pi_N$,

$$f(x) = \int_{M^{d-1}} f(y) K_{N,\eta}(\cos d(x, y)) \, d\sigma(y) = \sum_{j=1}^m w_j f(t_j) K_{N,\eta}(\cos d(x, t_j)).$$

This means that $f = T U_N(f)$ for $f \in \Pi_N$, where $U_N$ and $T$ are defined by (4.3) and (4.4), respectively. Thus, we can factor the identity $I : \Pi_N \cap L^p \to \Pi_{2N} \cap L^q$ as follows:

$$I : \quad \Pi_N \cap L^p \xrightarrow{U_N} \ell_p^m \xrightarrow{i_1} \ell_q^m \xrightarrow{i_2} \ell_{q,w}^m \xrightarrow{T} \Pi_{2N} \cap L^q,$$

where $\ell_{q,w}^m$ denotes the space $\mathbb{R}^m$ equipped with the norm $\|\cdot\|_{\ell_{q,w}^m}$ and $i_1, i_2$ both denote the identities. By well-known properties of $n$-widths (see [Pin, Chapter II]), it then follows that

$$S_n(\mathcal{B}_N^p, L^q) \leqslant S_n(\mathcal{B}_N^p, L^q \cap \Pi_{2N})$$

$$\leqslant \|U_N\|_{(\Pi_N \cap L^p, \ell_p^m)} \|i_2\|_{(\ell_q^m, \ell_{q,w}^m)} \|T\|_{(\ell_{q,w}^m, \Pi_{2N} \cap L^q)} S_n(b_p^m, \ell_q^m). \quad (4.6)$$

Since

$$0 \leqslant w_j \leqslant \frac{C_d}{N^{d-1}}, \quad 1 \leqslant j \leqslant m,$$

we have, for $x \in \mathbb{R}^m$,

$$\|x\|_{\ell_{q,w}^m} \leqslant C_d N^{-\frac{d-1}{q}} \|x\|_{\ell_q^m}. \quad (4.7)$$

Now a combination of (4.2), (4.5) and (4.7) gives

$$\|U_N\|_{(\Pi_N \cap L^p, \ell_p^m)} \|i_2\|_{(\ell_q^m, \ell_{q,w}^m)} \|T\|_{(\ell_{q,w}^m, \Pi_{2N} \cap L^q)} \leqslant C_d N^{(d-1)(\frac{1}{p}-\frac{1}{q})}.$$

The conclusion of Lemma 4.1 then follows from (4.6). □

**Proof of Theorem 1.1.** First, we prove the lower estimates. It will suffice to prove the lower estimates for $d_n(B_p^r, L^q)$, from which those of $\delta_n(B_p^r, L^q)$ will follow by the duality $\delta_n(B_p^r, L^q) = \delta_n(B_{q'}^r, L^{p'})$ and the inequality $\delta_n(B_p^r, L^q) \geqslant d_n(B_p^r, L^q)$. The lower estimates of $d_n(B_p^r, L^q)$ for $1 \leqslant p \leqslant 2 \leqslant q \leqslant \infty$ follow by the inequalities

$$d_n(B_p^r, L^q) \geqslant d_n(B_p^r, L^2) \geqslant C n^{-\frac{r}{d-1}+\frac{1}{p}-\frac{1}{2}} \tag{4.8}$$

and for $2 \leqslant p \leqslant q \leqslant \infty$ follow by

$$d_n(B_p^r, L^q) \geqslant d_n(B_p^r, L^p) \geqslant C n^{-\frac{r}{d-1}}, \tag{4.9}$$

where the last inequalities in (4.8) and (4.9) are contained in [BKLT, Theorem 1]. It remains to prove the lower estimates for $1 \leqslant q \leqslant p \leqslant \infty$. By the monotonicity of the norm $\|\cdot\|_q$, it will suffice to prove

$$d_n(B_\infty^r, L^1) \geqslant C n^{-\frac{r}{d-1}}, \tag{4.10}$$

with $C > 0$ being independent of $n$.

To prove (4.10), we assume $b_1^{-1} m^{d-1} \leqslant n \leqslant b_1 m^{d-1}$ with $b_1 > 0$ being independent of $n$ and $m$, and let $\{x_j\}_{j=1}^{2N_m} \subset M^{d-1}$ such that $N_m \asymp m^{d-1}$ and

$$B\left(x_i, \frac{2}{m}\right) \bigcap B\left(x_j, \frac{2}{m}\right) = \emptyset \quad \text{if } i \neq j.$$

We take $b_1 > 0$ sufficiently small so that $N_m \geqslant n$. Let $\varphi$ be a $C^\infty$-function on $\mathbb{R}$ supported in $[\frac{1}{2}, 1]$ and being equal to 1 on $[\frac{2}{3}, \frac{3}{4}]$. We define

$$\varphi_i(x) = \varphi(m \cdot d(x, x_i)), \quad 1 \leqslant i \leqslant 2N_m$$

and set

$$\mathcal{A}_m = \left\{ \sum_{j=1}^{2N_m} a_j \varphi_j(x) : \max_{1 \leqslant j \leqslant 2N_m} |a_j| \leqslant 1 \right\}.$$

Then, by the monotonicity of $d_n$, the relation

$$m^{-r}\mathcal{A}_m \subset CB_\infty^r \tag{4.11}$$

and the inequality

$$d_{N_m}(\mathcal{A}_m, L^1) \geqslant C, \tag{4.12}$$

with $C > 0$ depending only on $\varphi$ and $M^{d-1}$, will prove (4.10).

For the proof of (4.11), we note that in geodesic polar coordinates the Laplace–Beltrami operator $\triangle$ on $M^{d-1}$ equals $\triangle_\theta + \triangle'$, where $\triangle'$ denotes the Laplace–Beltrami operator on the sphere in $M^{d-1}$ of radius $\theta$ and

$$\triangle_\theta = \frac{\partial^2}{\partial\theta^2} + (p\cot\theta + 2q\cot 2\theta)\frac{\partial}{\partial\theta}$$

with $p, q$ depending only on $M^{d-1}$. See [He1, pp. 171–172]. So, by the definition, it follows that for $1 \leqslant i \leqslant 2N_m$

$$\|\triangle^v \varphi_i\|_\infty = \|\triangle_\theta^v \varphi_i\|_\infty \leqslant Cm^{2v}, \quad v = 1, 2, \ldots,$$

which, by a Kolmogorov type inequality (see [Di, Theorem 8.1]), implies

$$\|\triangle^{\frac{r}{2}} \varphi_i\|_\infty \leqslant C\|\varphi_i\|_\infty^{1-\frac{r}{2[\frac{r}{2}]+2}} \|\triangle^{1+[\frac{r}{2}]} \varphi_i\|_\infty^{\frac{r}{2[\frac{r}{2}]+2}} \leqslant Cm^r.$$

Since

$$\operatorname{supp}\varphi_i \bigcap \operatorname{supp}\varphi_j = \emptyset \quad \text{if } i \neq j, \tag{4.13}$$

(4.11) follows.

To prove (4.12), we note that for $1 \leqslant p \leqslant \infty$,

$$\|\varphi_1\|_p = \|\varphi_2\|_p = \cdots = \|\varphi_{2N_m}\|_p \asymp m^{-\frac{d-1}{p}} \asymp N_m^{-\frac{1}{p}}. \tag{4.14}$$

For $f \in L^1$ and $x \in M^{d-1}$, we define

$$P(f)(x) = \sum_{j=1}^{2N_m} \frac{\varphi_j(x)}{\|\varphi_j\|_2^2} \int_{M^{d-1}} f(y)\varphi_j(y)\,dy.$$

Then by (4.13) and (4.14), it follows that

$$\|P(f)\|_1 \leqslant C\|f\|_1, \qquad (4.15)$$

with $C > 0$ being an absolute constant. Let $X \subset L(M^{d-1})$ be an $N_m$-dimensional subspace. Using (4.13)–(4.15), we have

$$\sup_{f\in\mathcal{A}_m}\inf_{g\in X}\|f-g\|_1 \geqslant \sup_{f\in\mathcal{A}_m}\inf_{g\in X}\|f-P(g)\|_1 \geqslant CN_m^{-1}d_{N_m}(b_\infty^{2N_m}, \ell_1^{2N_m}).$$

This combined with the well-known inequality (see [LGM, p. 459])

$$d_{N_m}(b_\infty^{2N_m}, \ell_1^{2N_m}) \geqslant \frac{1}{2}N_m$$

gives (4.12), which completes the proof of the lower estimates.

Next, we prove the upper estimates. We shall prove the upper estimates only for the linear widths. The case of Kolmogorov widths can be treated similarly.

The upper estimates for $1 \leqslant q \leqslant p \leqslant \infty$ are obvious:

$$\delta_n(B_p^r, L^q) \leqslant \delta_n(B_p^r, L^p) \leqslant Cn^{-\frac{r}{d-1}}.$$

The last inequality is contained in [BKLT].

It remains to prove the upper estimates for $1 \leqslant p \leqslant 2 \leqslant q \leqslant \infty$. By the duality $\delta_n(B_p^r, L^q) = \delta_n(B_{q'}^r, L^{p'})$, it will suffice to prove them for $1 \leqslant p \leqslant 2 \leqslant q \leqslant p'$.

We define

$$A_0(f) = V_{1,\eta}(f), \quad A_j(f) = V_{2^j,\eta}(f) - V_{2^{j-1},\eta}(f) \quad \text{for } j \geqslant 1,$$

where the operators $V_{k,\eta}$ are defined by (2.11). Clearly, for $f \in W_p^r$,

$$f = \sum_{j=0}^\infty A_j(f)$$

and using Abel's transform finite times, we have

$$\|A_j(f)\|_p \leqslant C_{p,r} 2^{-jr} \left\|(-\Delta)^{\frac{r}{2}}(f)\right\|_p, \quad j \geqslant 0.$$

So by Lemma 4.1, it follows that for $1 \leqslant p \leqslant 2 \leqslant q \leqslant p'$,

$$\delta_n(B_p^r, L^q) \leqslant \delta_n(B_p^r, L^{p'}) \leqslant C_{p,r} \sum_{k=0}^{\infty} 2^{-kr} \delta_{n_k}(\mathcal{B}_{2^{k+1}}^p, L^{p'})$$

$$\leqslant C'_{p,r} \sum_{k=0}^{\infty} 2^{-kr} 2^{k(d-1)(\frac{2}{p}-1)} \delta_{n_k}(b_p^{m_k}, \ell_{p'}^{m_k}), \quad (4.16)$$

where $\sum_{k=0}^{\infty} n_k \leqslant n - 1$ and $m_k \asymp 2^{(d-1)k}$.

The rest of the proof is standard. See [Pin, pp. 236–241] and [Te, Chapter I, Section 4]. We sketch it as follows. Assume $C_1 2^{(d-1)v} \leqslant n \leqslant C_1^2 2^{(d-1)v}$ with $C_1 > 0$ to be specified later. We fix a real number $\rho \in \left(0, \frac{2r}{d-1} - \frac{2}{p}\right)$ and set

$$n_k = \begin{cases} m_k & \text{if } 0 \leqslant k \leqslant v, \\ \left[2^{(d-1)((1+\rho)v - k\rho)}\right] & \text{if } v < k < (1 + \rho^{-1})v, \\ 0 & \text{if } k \geqslant (1 + \rho^{-1})v. \end{cases}$$

Then a straightforward calculation shows $\sum_{k=0}^{\infty} n_k \asymp 2^{(d-1)v}$. Hence, we can take $C_1$ sufficiently large so that $\sum_{k=0}^{\infty} n_k \leqslant C_1 2^{(d-1)v} - 1 \leqslant n - 1$. It is clear that for $0 \leqslant k \leqslant v$,

$$\delta_{n_k}(b_p^{m_k}, \ell_{p'}^{m_k}) = \delta_{m_k}(b_p^{m_k}, \ell_{p'}^{m_k}) = 0$$

and for $k \geqslant (1 + \rho^{-1})v$,

$$\delta_{n_k}(b_p^{m_k}, \ell_{p'}^{m_k}) = \delta_0(b_p^{m_k}, \ell_{p'}^{m_k}) \leqslant 1.$$

For $v < k < (1 + \rho^{-1})v$, we use Gluskin's well-known estimates [Gl] to obtain

$$\delta_{n_k}(b_p^{m_k}, \ell_{p'}^{m_k}) \leqslant C_p m_k^{\frac{1}{p'}} n_k^{-\frac{1}{2}} \log^{\frac{1}{2}}\left(1 + \frac{m_k}{n_k}\right)$$

$$\leqslant C_{p,d} 2^{-\frac{(d-1)(1+\rho)v}{2}} 2^{(d-1)k(\frac{1}{p'}+\frac{\rho}{2})} \left((\rho+1)(k+1-v)\right)^{\frac{1}{2}}.$$

Now substituting these estimates into (4.16), we obtain, by straightforward calculation,

$$\delta_n(B_p^r, L^q) \leqslant C_p 2^{v\left(-r + (d-1)(\frac{1}{p} - \frac{1}{2})\right)} \leqslant C'_p n^{-\frac{r}{d-1} + \frac{1}{p} - \frac{1}{2}},$$

as required. This completes the proof. $\square$

## Acknowledgments

The authors would like to express their sincere gratitude to an anonymous referee for many helpful comments on this paper.